\documentclass[12pt,leqno]{article}
\usepackage{amsmath,amssymb,amscd}
\title{Presenting quantum Schur algebras as quotients of the quantized
universal enveloping algebra of ${\mathfrak {gl}}_2$}
\author{Stephen Doty and Anthony Giaquinto}
\date{December 1, 2000}

\newenvironment{pf}{{\em Proof.}}{\hfill$\square$\par\medskip}

\makeatletter %Redefine section and subsection formats
\renewcommand{\subsection}{\@startsection{subsection}{2}{0mm}{\baselineskip}{-\fontdimen2\font}{\normalfont\normalsize\bfseries\setcounter{equation}{0}}}
\makeatother

\newenvironment{xlem}{{\bf Lemma\ }\em}{\par\medskip}
\newenvironment{xthm}{{\bf Theorem\ }\em}{\par\medskip}

\newcommand{\N}{{\mathbb N}}
\newcommand{\Z}{{\mathbb Z}}
\newcommand{\Q}{{\mathbb Q}}
\newcommand{\A}{{\mathcal A}}
\newcommand{\B}{{\mathcal B}}
\newcommand{\End}{\operatorname{End}}
\renewcommand{\ker}{\operatorname{Ker}}

\newcommand{\gl}{{\mathfrak {gl}}}

\newcommand{\ts}{\otimes}
\newcommand{\Sv}{S_v(2,d)}
\newcommand{\Sz}{S_{\mathcal{A}}(2,d)}

\newcommand{\U}{\mathbf{U}}%{U_v(\mathfrak{gl}_2)}
\newcommand{\Uz}{\U_{\mathcal{A}}}%{U_{\mathcal{A}}(\mathfrak{gl}_2)}

\newcommand{\Ka}{{K_1}}
\newcommand{\Kb}{{K_2}}
\newcommand{\Kaa}{K_{1}^{-1}}
\newcommand{\Kbb}{K_{2}^{-1}}
\newcommand{\Ki}{{K_i}}
\newcommand{\KKa}{{\overline{K_1}}}
\newcommand{\KKb}{{\overline{K_2}}}
\newcommand{\vv}{v^{-1}}
\newcommand{\vbinom}[3]{\begin{bmatrix}#1; #2\\#3\end{bmatrix}}
\newcommand{\vbinomm}[2]{\begin{bmatrix}#1\\#2\end{bmatrix}}
\newcommand{\divided}[2]{#1^{(#2)}}
\newcommand{\K}{\begin{bmatrix}K_1\\b_1\end{bmatrix}\begin{bmatrix}K_2\\b_2\end{bmatrix}}
\newcommand{\KK}[2]{\begin{bmatrix}K_1\\#1 \end{bmatrix}\begin{bmatrix}K_2\\#2\end{bmatrix}}
\newcommand{\ea}{e^{(a)}}
\newcommand{\fc}{f^{(c)}}
\newcommand{\fa}{f^{(a)}}
\newcommand{\ec}{e^{(c)}}
\begin{document}
\maketitle
\begin{abstract} We obtain a presentation of the quantum Schur algebras
$\Sv$ by generators and relations.  This presentation is compatible
with the usual presentation of the quantized enveloping algebra $\U =
\mathbf{U}_v(\mathfrak{gl}_2)$.  In the process we find new bases for
$\Sv$.  We also locate the $\Z[v,v^{-1}]$-form of the quantum Schur
algebra within the presented algebra and show that it has a basis
which is closely related to Lusztig's basis of the $\Z[v,v^{-1}]$-form
of $\U$.

\end{abstract}

\section{Introduction}\label{AAA}
In \cite{DG} we gave a description of the rational Schur algebra
$S_{\Q}(2,d)$ in terms of generators and relations. This description
is compatible with the usual presentation of the universal enveloping
algebra $U(\mathfrak{gl}_2)$. We also described the integral Schur
algebra $S_{\Z}(2,d)$ as a certain subalgebra of the rational version
and an integral basis was exhibited. In this paper, we formulate and
prove quantum versions of those results.

Consider the Drinfeld-Jimbo quantized enveloping algebra $\U$
corresponding to the Lie algebra $\gl_2$. It has a natural
two-dimensional module $E$, and thus the tensor product $E^{\ts d}$ is
also a module for $\U$. Let
$$\rho_d:\U\rightarrow \End(E^{\ts d})$$ be the corresponding
representation. Amongst the various equivalent definitions of the
quantum Schur algebra $\Sv$ is that it is precisely the image of the
homomorphism $\rho_d$. (We shall not need to consider the other
definitions in this work.)  Since $\Sv$ is a homomorphic image of
$\U$, it is natural to ask for an efficient generating set of
$\ker(\rho_d)$, thereby giving a presentation of $\Sv$.

In what follows, we obtain a precise answer to this question.  Recall
that $\U$ is generated by elements $e$, $f$, $\Ka^{\pm 1}$, $\Kb^{\pm
1}$ subject to various well-known relations (see section 3). Now in
the representation $\rho_d$, it is easy to see that $\Ka\Kb=v^d$, so
we may use this relation to eliminate $\Kb$ (or $\Ka$) from the
generating set for $\Sv$.  Having done this, then our main result is
that the only additional relation needed to give the desired
presentation of $\Sv$ is the minimal polynomial of $\Ka$ in
$\End(E^{\ts d})$.

Set $\A=\Z[v,v^{-1}]$.  The algebra $\U$ contains a certain
$\A$-subalgebra $\Uz$ which may be viewed as a quantum version of the
integral form $U_{\Z}(\mathfrak{gl}_2)$ of the classical enveloping
algebra. The algebra $\Uz$, originally constructed by Lusztig, is
generated by the $v$-divided powers of $e$ and $f$ along with
$\Ka^{\pm 1}$ and $\Kb^{\pm 1}$. It has an $\A$-basis which is a
quantum analog of Kostant's basis for the integral form of the
classical enveloping algebra $U(\mathfrak{gl}_2)$.  The image of the
homomorphism $\rho_d$ upon restriction to $\Uz$ gives us an
``integral'' Schur algebra $\Sz$, which can be used to define a
version over any $\A$-algebra. We show that the integral Schur algebra
has a basis which is closely related to Lusztig's basis of $\Uz$.

Although the results of this paper are quantum versions of those
appearing in \cite{DG}, the techniques are somewhat different, since
some of the arguments given in that paper did not quantize
directly. In particular, the treatment here of the degree zero part
(generated by the images of $K_1^{\pm 1}$ and $K_2^{\pm 1}$) of $\Sv$
is totally different. Specifically, we exhibit an idempotent basis of
this subalgebra, whereas in \cite{DG} a PBW-type basis of the degree
zero part was used. The idempotent basis is more amenable to
computations and is precisely the kind of basis needed to handle the
general Schur algebras $S(n,d)$ and their quantizations. Another
difference between the results of this paper and the results
of \cite{DG} is that here we do not
obtain analogues of the ``restricted PBW-basis'', although we
believe such results should be true in the quantum case.

Finally, it is clear that analogous results hold in general, for any
$n$ and $d$, although one cannot expect to obtain such precise
reduction formulas in the general case as those given here and in
\cite{DG}.  The authors expect to treat the general case in a later
paper.

\section{Statement of results}\label{AA}
The main results of this paper are contained in the following
theorems. Let $\A=\Z[v,v^{-1}]$ with fraction field $\Q(v)$. Each
of the first three theorems gives a presentation of the quantum
Schur algebra $\Sv$ in terms of generators and relations. The
first result gives a presentation which is similar to that of
${\U}_v(\mathfrak{sl}_2)$.
\subsection{}\label{AAa0}
\begin{xthm}
Over $\Q(v)$, the quantum Schur algebra $\Sv$ is isomorphic to the
algebra generated by $e$, $f$, $K^{\pm 1}$ subject to the
relations:
\begin{align}
& KK^{-1}=K^{-1}K=1,\label{AAa0:a}\\
&K eK^{- 1} =v^{2}e,\qquad
K fK^{- 1} =v^{- 2}f,\label{AAa0:b}\\
& ef-fe=\frac{K - K^{-1}}{v-\vv},\label{AAa0:c}\\
& (K-v^d)(K-v^{d-2})\cdots (K-v^{-d+2})(K-v^{-d})=0.\label{AAa0:d}
\end{align}
\end{xthm}
The next two results give presentations of $\Sv$ which are
similar to that of ${\U}_v(\mathfrak{gl}_2)$.
\subsection{}\label{AAa}
\begin{xthm}
Over $\Q(v)$, the quantum Schur algebra $\Sv$ is isomorphic to the
algebra generated by $e$, $f$, $\Ka^{\pm 1}$ subject to the relations:
\begin{align}
& \Ka\Ka^{-1}=\Ka^{-1}\Ka=1,\label{AAa:a}\\
&\Ka e\Ka^{- 1} =ve,\qquad
\Ka f\Ka^{- 1} =v^{- 1}f\label{AAa:b}\\
& ef-fe=\frac{v^{-d}\Ka^2 - v^{d}\Ka^{-2}}{v-\vv},\label{AAa:c}\\
& (\Ka-1)(\Ka-v)(\Ka-v^2)\cdots (\Ka-v^d)=0.\label{AAa:d}
\end{align}
\end{xthm}
By a change of variable ($K_2=v^{d}\Kaa$) we obtain another
equivalent presentation of $\Sv$.
\subsection{}\label{AAb}
\begin{xthm}
Over $\Q(v)$, the quantum Schur algebra $\Sv$ is isomorphic to the
algebra generated by $e$, $f$, $\Kb^{\pm 1}$ subject to the relations:
\begin{align}
& \Kb\Kb^{-1}=\Kb^{-1}\Kb=1,\label{AAb:a}\\
&\Kb e\Kb^{- 1} =v^{- 1}e,\qquad
\Kb f\Kb^{- 1} =v f\label{AAb:b}\\
& ef-fe=\frac{v^d\Kb^{-2} - v^{-d}\Kb^{2}}{v-\vv},\label{AAb:c}\\
& (\Kb-1)(\Kb-v)(\Kb-v^2)\cdots (\Kb-v^d)=0.\label{AAb:d}
\end{align}
\end{xthm}
\noindent
For indeterminates $X,X^{-1}$ satisfying $X X^{-1} = X^{-1} X = 1$ and
any $t\in \N$ we formally set
$$\vbinomm{X}{t}=\prod_{s=1}^{t}
\frac{Xv^{-s+1}-X^{-1}v^{s-1}}{v^s-v^{-s}}.$$ This expression will
make sense if $X$ is replaced by any invertible element of a
$\Q(v)$-algebra.  The next result describes the $\A$-form $\Sz$ in
terms of the above generators, and gives an $\A$-basis for the
algebra.  In this we can take $\Sv$ to be given by either presentation
\ref{AAa} or \ref{AAb}, but we always assume that $\Ka$ and $\Kb$ are
related by the condition $\Ka\Kb=v^d$.

\subsection{}\begin{xthm}\label{AAc}
The integral Schur algebra $\Sz$ is isomorphic to the $\A$-subalgebra of
$\Sv$ generated by
$$e^{(m)}:=\frac{e^m}{[m]!},\quad f^{(m)}:=\frac{f^m}{[m]!}\quad
(m\in \N),\qquad \Ka^{\pm 1}.$$
The preceding statement is true when $\Ka$ is replaced by $\Kb$.
Moreover, an $\A$-basis for $\Sz$ is the set consisting of all
$$ \ea \K \fc $$
such that the natural numbers $a,b_1,b_2,c$ are constrained by
the conditions $a+b_1+c\leq d,\ b_1+b_2=d$.
Another such basis consists of all
$$f^{(a)} \K e^{(c)} $$
such that the natural numbers $a,b_1,b_2,c$ satisfy the constraints
$a+b_2+c\leq d,\  b_1+b_2=d$.
\end{xthm}

The following reduction formulas make it possible, in principle, to
express the product of two basis elements as an $\A$-linear
combination of basis elements.

\subsection{}\begin{xthm}\label{AAd}
In $\Sz$ we have the following reduction formulas for all $s\ge1$ and all
$a,b_1,b_2,c\in \N$ with $b_1+b_2=d$:
\begin{gather}
\begin{aligned}\ea\K\fc=& \\
\sum_{k=s}^{\min(a,c)}(-1)^{k-s}&\vbinomm{k-1}{s-1}\vbinomm{b_1+k}{k}
\divided{e}{a-k}\KK{b_1+k}{b_2-k}\divided{f}{c-k}
\end{aligned}\label{AAd:a}\\
\begin{aligned}
f^{(a)}\K e^{(c)}=& \\
\sum_{k=s}^{\min(a,c)}(-1)^{k-s}&\vbinomm{k-1}{s-1}\vbinomm{b_2+k}{k}
\divided{f}{a-k}\KK{b_1-k}{b_2+k}\divided{e}{c-k}
\end{aligned}\label{AAd:b}
\end{gather}
where $s=a+b_1+c-d$ in (a) and $s=a+b_2+c-d$ in (b).
\end{xthm}

The next theorem, the quantum analogue of \cite[Thm.\ 2.4]{DG},
provides yet another kind of basis for $\Sz$.  In this case we will
deduce it as a direct consequence of Theorem \ref{AAc}, while in
\cite{DG} the analogue of the idempotent basis in Theorem \ref{AAc}
was not needed.  We remark that the analogues of the integral
idempotent bases in Theorem \ref{AAc} above do hold in the classical
situation.
\subsection{}\begin{xthm}\label{AAe}
The set $$\left\{\ea \vbinomm{K_1}{b} \fc \quad\mid\quad a,b,c\in \N,\
a+b+c\leq d\right\}$$ is an $\A$-basis for $\Sz$. Another such basis
is given by the set $$\left\{\divided{f}{a} \vbinomm{K_2}{b}
\divided{e}{c} \quad\mid\quad a,b,c\in \N,\ a+b+c\leq d\right\}.$$
\end{xthm}

\section{Quantized enveloping algebras}\label{A}
\subsection{}\label{Aa}
The Drinfeld-Jimbo quantized enveloping algebra
$\U=\U_v(\mathfrak{gl}_2)$ is defined to be the $\Q(v)$-algebra with
generators $e$, $f$, $\Ka$, $\Kaa$, $\Kb$, $\Kbb$ and relations
\begin{align}
& \Ka\Kb=\Kb\Ka,\label{Aa:a}\\
& \Ki\Ki^{-1}=\Ki^{-1}\Ki=1 \quad (i=1,2),\label{Aa:b}\\
& \Ka e\Kaa =ve, \quad \Ka f\Kaa=v^{-1}f,\label{Aa:c}\\
& \Kb e\Kbb =\vv e, \quad \Kb f\Kaa=vf,\label{Aa:d}\\
& ef-fe=\frac{\Ka\Kbb - \Kaa\Kb}{v-\vv}\label{Aa:e}.
\end{align}
Let $\U^+$ (respectively $\U^-$) be the subalgebra generated by $e$
(respectively $f$) and let $\U^0$ be the subalgebra generated by
$K_1^{\pm 1}$, $K_2^{\pm 1}.$ There are $\Q(v)$-vector space
isomorphisms
$$\U \cong \U^+\ts \U^0 \ts \U^- \cong \U^-\ts \U^0 \ts \U^+ $$ and
moreover it is well-known that the algebra $\U$ has PBW-type bases
$\left\{ e^a \Ka^{b_{1}} \Kb^{b_{2}}f^c\right\}$ and $\left\{ f^a
\Ka^{b_{1}} \Kb^{b_{2}}e^c\right\}$ where $a,c\in \N$ and $b_1,b_2\in
\Z.$

In the algebra $\U$, set $K=K_1K_2^{-1}$, and define
$\U_v(\mathfrak{sl}_2)$ to be the
subalgebra of $\U$
generated by $e$, $f$, and $K^{\pm 1}$. The familiar relations
that these elements satisfy are easily deducible from \eqref{Aa:a}
- \eqref{Aa:e}.

\subsection{}\label{Ab}
For $r,s \in \Z$ define
\begin{align}
& [r]=\frac{v^r-v^{-r}}{v-\vv}\label{Ab:a}\\
& \vbinomm {r}{s} = \frac{[r][r-1]\cdots [r-s+1]}{[1][2]\cdots
[s]}.\label{Ab:b}\end{align}
These satisfy the well-known
identities
\begin{align}
& [r+s]=v^{-s}[r]+v^r [s]\label{Ab:c}\\
& \vbinomm {r+1}{s}=v^{-s}\vbinomm{r}{s}+v^{r-s+1}\vbinomm{r}{s-1}.
 \end{align}
For every $m,t\in \N,$ $c\in \Z$ and any
element $X$ an $\Q(v)$-algebra
define
 \begin{align}
 & [m]!=[m][m-1]\cdots [1]\label{Ab:f}\\
 & X^{(m)}=\frac{X^m}{[m]!}\label{Ab:g}\\
 & \vbinom {X}{c}{t} = \begin{cases} \displaystyle{ \prod_{i=1}^t
\frac{Xv^{c-i+1}-X^{-1}v^{-c+i-1}}{v^i-v^{-i}} }
 \quad \text{if}\quad t\neq 0 \\
1 \quad \text{if} \quad t=0.\end{cases}\quad
(X\quad \text{invertible}) \label{Ab:h}
\end{align}
\noindent We note that $\vbinomm{X}{t}$, as defined in section 2,
coincides with the element $\vbinom{X}{0}{t}$.

In \cite{Lu}, Lusztig investigates many relations which hold
in quantized enveloping algebras. Those which we will need
are contained in the following Lemma.
\subsection{} \begin{xlem}\label{Ad}
For any $c,n\in \Z$ and $m,t,t'\in \N$ the following identities hold
in $\U$:
\begin{gather}
\Ka^n\,e\, K_1^{-n}=v^n\,e, \quad \Kb^n\,e\, K_2^{-n}=v^{-n}\,e, \label{Ad:a}\\
\Ka^n\,f\, K_1^{-n}=v^{-n}\,f, \quad \Kb^n\,f\, K_2^{-n}=v^{n}\,f,\label{Ad:b}\\
\vbinom{\Ka}{c}{t} \,e=e\, \vbinom{\Ka}{c+1}{t},\quad
\vbinom{\Ka}{c}{t} \,f=f\, \vbinom{\Ka}{c-1}{t},\label{Ad:c}\\
\vbinom{\Kb}{c}{t} \,e=e\, \vbinom{\Kb}{c-1}{t},\quad
\vbinom{\Kb}{c}{t} \,f=f\, \vbinom{\Kb}{c+1}{t},\label{Ad:d}\\
\divided{f}{m}\, e= e\, \divided{f}{m} - \vbinom{\Ka\Kbb}{m-1}{1}\divided{f}{m-1},\label{Ad:e}\\
f\, \divided{e}{m}= \divided{e}{m}\, f -\divided{e}{m-1}\vbinom{\Ka\Kbb}{m-1}{1},\label{Ad:f}\\
\vbinom{\Ki}{c+1}{t+1}=v^{t+1}\vbinom{\Ki}{c}{t+1}
+v^{t-c}\Ki ^{-1}\vbinom{\Ki}{c}{t},\label{Ad:g}\\
\vbinomm{\Ki}{t}\vbinom{\Ki}{-t}{t'}=\vbinomm{t+t'}{t}\vbinomm{\Ki}{t+t'},\label{Ad:h}\\
\vbinom{\Ki}{c}{t}=\sum_{j=0}^t v^{c(t-j)}\vbinomm{c}{j}\Ki
^{-j}\vbinomm{\Ki}{t-j}\quad (c\geq 0).\label{Ad:i}
\end{gather}
\end{xlem}

\begin{pf}
All of these identities are special cases of those appearing on
pages 269-270 of \cite{Lu}.
\end{pf}

\section{The algebra $\B_d$}
In this section we define a homomorphic image of $\U$ which will
turn out to be isomorphic to the quantum Schur algebra $\Sv$.

\subsection{}\label{Ba}
Let $d$ be a fixed nonnegative integer and define $\B_d$ to be the
$\Q(v)$-algebra generated by $e$, $f$, $\Ka ^{\pm 1}$,
$\Kb ^{\pm 1}$ subject to relations
\ref{Aa}\eqref{Aa:a}-\eqref{Aa:e}, along with the additional
relations
\begin{align}
&\Ka \Kb = v^d \label{Ba:a}\\
& (\Ka-1)(\Ka-v)\cdots (\Ka-v^d)=0\label{Ba:b}
\end{align}
Note that relation (b)
can be replaced by the equivalent relation
\begin{equation}\label{Ba:c}
(\Kaa-1)(\Kaa-\vv)\cdots (\Kaa-v^{-d})=0
\end{equation}
and in the presence of relation (a) it can also be replaced by either
\begin{align}
 & (\Kb-1)(\Kb-v)\cdots (\Kb-v^{d})=0,\quad {\text {or}}\label{Ba:d}\\
 & (\Kbb-1)(\Kbb-\vv)\cdots (\Kbb-v^{-d})=0.\label{Ba:e}
\end{align}

The defining relations of $\B_d$ are invariant if $e$ and $f$ are
interchanged along with $K_1$ and $K_2$. These interchanges therefore
induce an automorphism of $\B_d$. We shall often make use of this
property, which we will call {\it symmetry}, in the sequel.

Let $\B_d^0$ be the subalgebra of $\B_d$ generated by $\Ka ^{\pm 1}$,
$\Kb ^{\pm 1}.$ It follows from relations
\ref{Ba}\eqref{Ba:a}-\eqref{Ba:b} that $\dim(\B_d^0)=d+1$.
\subsection{} \begin{xlem}\label{Bb}
In the algebra $\B_d^0$ we have
$$\vbinomm{\Ka}{d+1}=0=\vbinomm{\Kb}{d+1},\qquad
\vbinomm{\Kaa}{d+1}=0=\vbinomm{\Kbb}{d+1}.$$
\end{xlem}

\begin{pf} Using definition \ref{Ab}\eqref{Ab:h} we have that
\begin{align*}
\vbinomm{\Ka}{d+1}&=\prod_{i=1}^{d+1}
\frac{\Ka v^{1-i}-\Kaa v^{i-1}}{v^i-v^{-i}}\\
&=\prod_{i=1}^{d+1} \frac{(\Kaa
v^{1-i})(K_1^2-v^{2(i-1)})}{v^i-v^{-i}}\\
&=\prod_{i=1}^{d+1} \frac{(\Kaa
v^{1-i})(K_1-v^{(i-1)})(K_1+v^{(i-1)})}{v^i-v^{-i}}\\
&=0\quad {\text{by relation \ref{Ba}\eqref{Ba:b}.}}
\end{align*}
The other equalities follow in a similar manner using
\ref{Ba}\eqref{Ba:c}-\eqref{Ba:e}.
\end{pf}
More generally we have
\subsection{} \begin{xlem}\label{Bc}
In the algebra $\B_d^0$ we have
$$\vbinomm{\Ka}{b_1}\vbinomm{\Kb}{b_2}=0$$
whenever $b_1+b_2=d+1.$ \end{xlem}

\begin{pf}
\begin{align*}
\vbinomm{\Ka}{b_1}\vbinomm{\Kb}{b_2}&=\prod_{i=1}^{b_1}
\frac{\Ka v^{1-i}-\Kaa v^{i-1}}{v^i-v^{-i}}
\prod_{i=1}^{b_2}
\frac{\Kb v^{1-i}-\Kbb v^{i-1}}{v^i-v^{-i}}\\
&=\prod_{i=1}^{b_1}
\frac{\Ka v^{1-i}-\Kaa v^{i-1}}{v^i-v^{-i}}
\prod_{i=1}^{b_2}
\frac{\Kaa v^{d+1-i}-\Kb v^{-d+i-1}}{v^i-v^{-i}}\\
&=(-1)^{b_2}\vbinomm{d+1}{b_2}\vbinomm{K_1}{d+1}\\
&=0 \quad {\text{by
\ref{Ba}\eqref{Ba:b}}}.\end{align*}
\end{pf}
An immediate consequence of the preceding lemma is that
$$\vbinomm{\Ka}{b_1}\vbinomm{\Kb}{b_2}=0$$
whenever $b_1+b_2\ge d+1.$ When $b_1+b_2\leq d$, these elements
are non-zero, but we shall see that the most important case is
$b_1+b_2=d.$

\subsection{} \begin{xlem}\label{Bd}
Suppose that $b_1+b_2=d$ and $t\in \N.$ Then for $i=1,2$ the following
identities hold in $\B_d^0$:
\begin{gather}
K_i \K=v^{b_i}\K \label{Bd:a}\\
\vbinom{K_i}{c}{t} \K=\vbinomm{b_i+c}{t}\K \label{Bd:b}
\end{gather} \end{xlem}

\begin{pf}
The second equality of the lemma follows immediately
from the first using definition \ref{Ab}\eqref{Ab:h}, and so we only need
to prove (a). Consider the case $i=1$. We have
\begin{align*}
(\Ka &- v^{b_1})\K \\
& = (\Ka-v^{b_1})  \prod_{i=1}^{b_1}
\frac{\Ka v^{1-i}-\Kaa v^{i-1}}{v^i-v^{-i}}
 \prod_{i=1}^{b_2}
\frac{\Kb v^{1-i}-\Kbb v^{i-1}}{v^i-v^{-i}}\\
& = (-1)^{b_2}(\Ka -v^{b_1}) \prod_{\substack{ i=1\\i\neq
b_1+1}}^{d+1}
\frac{\Ka v^{1-i}-\Kaa v^{i-1}}{v^i-v^{-i}}\\
&= (-1)^{b_2}\,(\Ka-v^{b_1})\prod_{\substack{ i=1\\i\neq
b_1+1}}^{d+1}
 \frac{v^{1-i}\Kaa (\Ka^2-v^{2(i-1)})}{v^i-v^{-i}}\\
&=0\quad {\text {by \ref{Ba}\eqref{Ba:b}.}}
\end{align*}
This proves identity (a) for $i=1.$ The case $i=2$ follows
from symmetry.
\end{pf}
The following result establishes the structure of the algebra
$\B_d^0$; it will later be crucial in determining the structure of
the entire algebra $\B_d$.
\subsection{} \begin{xthm}\label{Be}
In the algebra $\B_d^0$, the set
$$\left\{ \K
\quad | \quad b_1+b_2=d\right\}$$
is a basis of mutually orthogonal idempotents whose sum is the
identity.
\end{xthm}
\begin{pf}
By Lemma \ref{Bd} we have that
$$\left(\K \right)^2=\vbinomm{b_1}{b_1}\vbinomm{b_2}{b_2}\K=\K$$
and so each $\K$ is an idempotent.
Now suppose that $b_1+b_2=b_1'+b_2'=d$ and
$b_1\neq b_1'.$ Then either $b_1+b_2'\geq d+1$ or $b_1'+b_2\geq
d+1$ and so orthogonality follows by Lemma \ref{Bc}.
Thus in the algebra $\B_d^0$ (which has dimension $d+1$), we have a
set of $d+1$ distinct mutually orthogonal idempotents. It follows
that these must form a basis and that their sum is the identity.
\end{pf}

These idempotents have pleasant commutation relations with the
elements $e$ and $f$, given in the next lemma.

\subsection{}\begin{xlem}\label{Bf}
Suppose that $b_1+b_2=d$ and $a\in \N.$ Then in the algebra $\B_d$
\begin{align}
& \K e^a =\begin{cases}
e^a\,\vbinomm{\Ka}{b_1-a}\vbinomm{\Kb}{b_2+a}\quad \text{if}\quad
b_1\geq a\\ \\
0\quad \text{if}\quad b_1<a \end{cases}\\
& e^a \K  =\begin{cases}
\vbinomm{\Ka}{b_1+a}\vbinomm{\Kb}{b_2-a} e^a \quad \text{if}\quad
b_2\leq a\\ \\
0\quad \text{if}\quad b_2>a \end{cases}\\
& f^a \K =\begin{cases}
\vbinomm{\Ka}{b_1-a}\vbinomm{\Kb}{b_2+a}f^a\quad \text{if}\quad
b_1\geq a\\ \\
0\quad \text{if}\quad b_1<a.  \end{cases}\\
& \K f^a =\begin{cases}
f^a \vbinomm{\Ka}{b_1+a}\vbinomm{\Kb}{b_2-a}\quad \text{if}\quad
b_2\leq a\\ \\
0\quad \text{if}\quad b_2>a.  \end{cases}
\end{align}
\end{xlem}
\begin{pf}
Each of the relations is similar to prove, so we will only verify
the first one.
By \ref{Ad}\eqref{Ad:c} and \ref{Ad}\eqref{Ad:d} we have $$\K \,
e^a=e^a\, \vbinom{\Ka}{a}{b_1} \vbinom{\Kb}{-a}{b_2}$$ and thus
$$\vbinomm{b_2+a}{a}^{-1}\K \,e^a\,\vbinomm{\Kb}{a} =
\vbinomm{b_2+a}{a}^{-1} e^a\, \vbinom{\Ka}{a}{b_1}
\vbinom{\Kb}{-a}{b_2}\vbinomm{\Kb}{a}.$$ Now using
\ref{Ad}\eqref{Ad:d} and  \ref{Ad}\eqref{Ad:h}, this equality
becomes $$\vbinomm{b_2+a}{a}^{-1}\K \,\vbinom{\Kb}{a}{a}\, e^a=
\vbinomm{b_2+a}{a}^{-1}\,e^a\, \vbinom{\Ka}{a}{b_1}\,
\vbinomm{\Kb}{b_2+a}.$$
Transforming this further using \ref{Ad}\eqref{Ad:i} and
\ref{Bd}\eqref{Bd:b} yields
$$ \K \, e^a =\vbinomm{b_2+a}{a}^{-1}\, e^a \left( \sum_{j=0}^{b_1}
v^{a(b_1-j)}\vbinomm{a}{j}\Ka ^{-j}\vbinomm{\Ka}{b_1-j}\right)
\vbinomm{\Kb}{b_2+a}.$$
Since $\vbinomm{\Ka}{b_1-j}\vbinomm{\Kb}{b_2+a}=0$ whenever $j<a$
and $\vbinomm{a}{j}=0$ whenever $j>a$, the only possible non-zero
term on the right side of the equality corresponds to the case
$j=a.$ Consequently, $\K\, e^a = 0$ whenever $b_1<a$ and
if $b_1\geq a$ then
$$\K \, e^a=v^{a(b_1-a)}\Ka
^{-a}\vbinomm{\Ka}{b_1-a}\vbinomm{\Kb}{b_2+a}.$$ The remaining part of
claim (a) of the lemma now follows using \ref{Bd}\eqref{Bd:a}.
\end{pf}

\subsection{Remark} An immediate consequence of the preceding lemma is
that
$$e^a\K =0=\K f^a\quad  \text{whenever}\quad  a+b_1\geq d+1$$
and
$$f^a\K =0=\K e^a\quad \text{whenever}\quad a+b_2\geq d+1.$$
In particular, this implies that both $e$ and $f$ are nilpotent
of index $d+1$.

\medskip

The results obtained thus far show that the sets
$$\left\{ \ea\, \K \, \fc \mid a,c\leq d, b_1+b_2=d \right\}$$
and
$$\left\{ f^{(a)}\, \K \, \fc \mid a,c\leq d, b_1+b_2=d \right\}$$ are
spanning sets for the algebra $\B_d.$ We will henceforth refer to
elements of these spanning sets simply as {\em monomials}. These sets,
however, are not bases. To establish this fact, we need some
terminology. Define the {\it fake degree} of a monomial $\ea\, \K \,
\fc$ (respectively $\fa\, \K \, \ec$) to be $a+b_1+c$ (respectively
$a+b_2+c$) and, in both cases, define its {\it height} to be $a+c$.
\subsection{}\begin{xthm}\label{Bg}
Suppose that $b_1+b_2=d$. Then in the algebra $\B_d$ all monomials of the form
$\ea\, \K \, \fc$ (respectively $\fa\, \K \, \ec$)
of fake degree $d+1$ are expressible as $\Q(v)$-linear combinations of
monomials of the same form of strictly smaller fake degree and height.
\end{xthm}

\begin{pf}
By symmetry it is enough to prove the claim for monomials of the
form $\ea\, \K \, \fc$.
We use induction on height. The base case here is height one since
there are no height zero monomials of fake degree $d+1$. Consider
a monomial of fake degree $d+1$ and height $s\geq 1$. Suppose that
$a\geq c$. Then $a\geq 1$ and by Lemma \ref{Bf} we have
$$\divided{e}{a-1}\KK{b_1+1}{b_2-1} \,e\,\fc=[a]\ea \K \fc$$
and so the claim will follow once we show that
$\divided{e}{a-1}\KK{b_1+1}{b_2-1} \,e\,\fc$ can be expressed in
the desired form. But, by \ref{Ad}\eqref{Ad:e} and
\ref{Bd}\eqref{Bd:b} we have
\begin{align*}
\divided{e}{a-1}&\KK{b_1+1}{b_2-1} \,e\,\fc \\
&=\divided{e}{a-1}\KK{b_1+1}{b_2-1}\left( \fc e +\vbinom{\Ka
\Kbb}{c-1}{1} \divided{f}{c-1}\right)\\
&= \divided{e}{a-1}\KK{b_1+1}{b_2-1} \fc e \\
& \qquad    + [b_1-b_2+c-1]\divided{e}{a-1}\KK{b_1+1}{b_2-1}\divided{f}{c-1}
\end{align*}
where the second term on the right hand side is zero if $c=0.$
However, in general, the second term is a monomial of fake degree $d$
and height $s-2$, which is of the desired form. To analyze the first
term, we need to use the fact that if $M$ is a monomial of fake degree
$d'$ and height $s'$, then $Me$ is an $\Q(v)$-linear combination of
monomials of fake degree at most $d'$ and height at most $s'+1$ (this
claim can be verified using \ref{Ad}\eqref{Ad:e} and
\ref{Bd}\eqref{Bd:b}). Now by induction,
$\divided{e}{a-1}\KK{b_1+1}{b_2-1}\fc$ is expressible as an
$\Q(v)$-linear combination of monomials of fake degree at most $d$ and
height at most $s-2$. Thus, by our claim,
$\divided{e}{a-1}\KK{b_1+1}{b_2-1}\fc e$ is expressible as a linear
combination of terms of fake degree at most $d$ and height at most
$s-1$. This completes the proof in the case $a\geq c$; the case $a\leq
c$ is similar and is omitted.
\end{pf}

\section{Identifications} In this section we show that $\B_d$ is
isomorphic to $\Sv$.

\subsection{}\label{Ae}
Let $E$ be an $\Q(v)$-module with basis $\{e_1,e_2\}$. There is a
canonical representation $\rho: \U\rightarrow \End_{\Q(v)}(E)$
defined by:
$$e \mapsto \begin{pmatrix} 0&1\\ 0&0\end{pmatrix},\quad
f \mapsto \begin{pmatrix} 0&0\\ 1&0\end{pmatrix},\quad
K_1 \mapsto \begin{pmatrix} v&0\\ 0&1\end{pmatrix},\quad
K_2 \mapsto \begin{pmatrix} 1&0\\ 0&v\end{pmatrix}.$$
Since $\U$ is a bialgebra, we obtain a representation
$$\rho_d: \U\rightarrow \End_{\Q(v)}(E^{\ts d})$$ in the $d$th tensor
power of $E$, for every $d\in \N$.  Specifically, $\rho_d =\rho^d
\circ \Delta ^{d-1}$, where $\Delta^{d-1}:\U\rightarrow \U^{\ts d}$ is
the iterated comultiplication map. As stated earlier, $\Sv$ is the
image of the homomorphism $\rho_d.$

\subsection{}\begin{xlem}\label{Ca}
Write $\overline{X}$ for the image of $X\in \U$ under the
representation $\rho_d: \U\rightarrow \End_{\Q(v)}(E^{\ts d})$. Then
we have the identities
\begin{gather}
\KKa \KKb = v^d \\
(\KKa-1)(\KKa-v)\cdots (\KKa-v^d)=0
\end{gather}
\end{xlem}
\begin{pf}
We have $\KKa \KKb = \rho_d(\Ka \Kb)=\rho_1(\Ka \Kb)^{\ts d}$
since $\Delta (\Ki)=\Ki \ts \Ki$. But by \ref{Ae},
$$\rho_1(\Ka \Kb)= \begin{pmatrix} v&0\\ 0&1\end{pmatrix}
\begin{pmatrix} 1&0\\ 0&v\end{pmatrix}=v\begin{pmatrix} 1&0\\ 0&1\end{pmatrix}
=v1_V$$
and so $\rho_1(\Ka \Kb)^{\ts d}=v^d1_{V^{\ts d}}$, which proves
claim (a). The case $d=1$ for claim (b) follows immediately since
$\rho_1(\Ka)=\begin{pmatrix} v&0\\ 0&1\end{pmatrix}$.  Now $\rho_d
(\Ka)=\Ka ^{\ts d}$ is a diagonal matrix whose entries are
$1,v,\ldots, v^d$ (not counting multiplicities). The claim
follows. \end{pf}
\noindent Henceforth we shall omit the bar, writing the image of $X$
simply as $X$ instead of $\overline{X}$.

\subsection{}\begin{xthm}\label{Cb}
The Schur algebra $\Sv$ is isomorphic as a $\Q(v)$-algebra to
$\B_d$. Moreover, the set of all monomials $\ea \K \fc$ ($a+b_1+c\leq
d$) is a basis over $\Q(v)$. Similarly, another such basis consists of
all the monomials $\fa \K \ec$ ($a+b_2+c\leq d$). \end{xthm}
\begin{pf}
By the preceding lemma we see that the surjection $\rho_d:\U
\to \Sv$ factors through $\B_d$, giving the commutative diagram
\begin{equation}
\begin{CD}
\U @>>> \Sv\\
@VVV           @AAA\\
\B_d @= \B_d
\end{CD}
\end{equation}
in which all arrows are surjections. By Theorem \ref{Bg}, the set
of all monomials of the form $\ea \K \fc$ satisfying $a+b_1+c\leq
d$ spans the algebra $\B_d$. This spanning set is in one-to-one
correspondence with the set of all monomials in 4 commuting
variables of total degree $d$ (set one of the variables equal to 1
to get the correspondence). It is well known that this number is
the dimension of $\Sv$. Hence the map $\B_d \rightarrow \Sv$ must
be an isomorphism and the spanning set is a basis. A similar
argument establishes that the set of all $\fa\, \K \, \ec$
satisfying $a+b_2+c\leq d$ is also a basis of $\B_d$. This
completes the proof. 
\end{pf}

\subsection{} Theorems \ref{AAa} and \ref{AAb} follow immediately
from the previous result by using the relation $\Ka\Kb=v^d$ to remove
either $\Ka$ or $\Kb$ from the generating set. However, the results of
the preceding theorem are not sufficient to establish Theorem
\ref{AAc} since we only know that the coefficients in the reduction
formulas are elements of $\Q(v)$, not necessarily of $\A$.  In the
next section we indeed show that these coefficients all lie in
$\A$. The process by which we prove this does not rely on Theorem
\ref{Cb} and so the forthcoming results give alternative proofs of
Theorems \ref{AAa} and \ref{AAb} as well as a proof of \ref{AAc}.

We now prove Theorem \ref{AAa0}. In the algebra $\Sv$, set
$K=v^{-d}K_1^2$. Using Lemma \ref{Bd} and Theorem \ref{Be} we obtain
the equality
$$K_1^2=v^dK=\sum
_{b_1+b_2=d}v^{2b_1}\vbinomm{\Ka}{b_1}\vbinomm{\Kb}{b_2} =
\left(\sum
_{b_1+b_2=d}v^{b_1}\vbinomm{\Ka}{b_1}\vbinomm{\Kb}{b_2}\right)^2.
$$
Thus, $e$, $f$, and $K^{\pm 1}$ comprise
another generating set for $\Sv$. Relations (a)-(c) of
Theorem \ref{AAa0} are easily seen to be equivalent to the
corresponding parts of Theorem \ref{AAa}. Relation (d) of these
theorems are also equivalent. To see this, note that
\begin{align*}
\prod _{i=0}^d(K-v^{d-2i})&=v^{-d^2-d}\prod _{i=0}^d(K_1^2-v^{2(d-i)})\\
&=v^{-d^2-d}\prod _{i=0}^d(K_1-v^{(d-i)})(K_1+v^{(d-i)}).
\end{align*}
But by Theorem \ref{AAa} this expression is equal to zero.
This completes the proof of Theorem \ref{AAa0}.

\renewcommand{\K}{K_{b_1,b_2}}
\renewcommand{\KK}[2]{K_{#1,#2}}
\section{The integral reduction}
In this section we prove Theorems \ref{AAc}, \ref{AAd}, and \ref{AAe}.
For typeographic convenience, we will throughout this section use the
abbreviation $\K$ in place of the more cumbersome notation
${\begin{bmatrix}K_1\\b_1\end{bmatrix}\begin{bmatrix}K_2\\b_2\end{bmatrix}}$.

\subsection{}\begin{xthm}\label{Da}
In the algebra $\B_d$ we have the equality $$\divided{e}{a}\K
\divided{f}{c}=\sum_{k=1}^{\min(a,c)}(-1)^{k-1}
\vbinomm{b_1+k}{k}\divided{e}{a-k}\KK{b_1+k}{b_2-k}\divided{f}{c-k}$$
for all $a,b_1,b_2,c\in \N$ satisfying $a+b_1+c=d+1$ and $b_1+b_2=d$.
Similarly we have the equality $$f^{(a)}\K e^{(c)}=
\sum_{k=1}^{\min(a,c)}(-1)^{k-1}\vbinomm{b_2+k}{k}
\divided{f}{a-k}\KK{b_1-k}{b_2+k}\divided{e}{c-k}$$ for all
$a,b_1,b_2,c\in \N$ satisfying $a+b_2+c=d+1$ and $b_1+b_2=d$.

\end{xthm}
\begin{pf}
By symmetry it suffices to prove the first reduction formula only.
Suppose first that $c=0$. Then the claim follows from the remark
following Lemma \ref{Bf}. If $c>0$ we will prove the desired
result by induction on $a$. The base case here is $a=0$, which
also follows from the remark mentioned above. The identity to be
proved can be rewritten in the form
\begin{equation}\label{Da:a}
0=\sum_{k=0}^{\min(a,c)}(-1)^k
\vbinomm{b_1+k}{k}\divided{e}{a-k}\KK{b_1+k}{b_2-k}\divided{f}{c-k}.
\end{equation}
Suppose this equation is satisfied for some fixed quadruple
$a,b_1,b_2,c$ satisfying the conditions $a+b_1+c=d+1$, $b_1+b_2=d$
and $c\geq 0$. From this we will derive the result for the case
$a+1, b_1-1,b_2+1,c$. The idea is to multiply \eqref{Da:a} on the
right by $e$ and commute it all the way to the left. Using Lemma
\ref{Ad}\eqref{Ad:e} we obtain
\begin{align*}
0&=\sum_{k=0}^{\min(a,c)}(-1)^k
\vbinomm{b_1+k}{k}\divided{e}{a-k}\KK{b_1+k}{b_2-k}\divided{f}{c-k}e\\
&=\sum_{k=0}^{\min(a,c)}(-1)^k
\vbinomm{b_1+k}{k}\divided{e}{a-k}\KK{b_1+k}{b_2-k}\\
&\hspace{2in} \times \left(
e\divided{f}{c-k}-\vbinom{\Ka\Kbb}{c-k-1}{1}\divided{f}{c-k-1}\right).
\end{align*}
Then by Lemmas \ref{Bd} and \ref{Bf} this equation becomes
\begin{align*}
0&=\sum_{k=0}^{\min(a,c)}(-1)^k
\vbinomm{b_1+k}{k}\divided{e}{a-k}e\KK{b_1+k-1}{b_2-k+1}\divided{f}{c-k}\\
&-\sum_{k=0}^{\min(a,c)}(-1)^k
\vbinomm{b_1+k}{k}[b_1-b_2+k+c-1]\divided{e}{a-k}\KK{b_1+k}{b_2-k}
\divided{f}{c-k-1}
\end{align*}
which is equivalent to the identity
\begin{align*}
0&=\sum_{k=0}^{\min(a,c)}(-1)^k
\vbinomm{b_1+k}{k}[a-k+1]\divided{e}{a-k+1}\KK{b_1+k-1}{b_2-k+1}
\divided{f}{c-k}\\
&+\sum_{k=1}^{1+\min(a,c)}(-1)^{k-1}
\vbinomm{b_1+k-1}{k-1}[a-b_1-k+1]\\
&\hspace{2in}\times \divided{e}{a-k+1}\KK{b_1+k-1}{b_2-k+1}
\divided{f}{c-k}.
\end{align*}
(The equality $b_1-b_2+k+c-1=-(a-b_1-k)$ follows from the
hypotheses $a+b_1+c=d+1$ and $b_1+b_2=d$.) Writing $m$ for
$\min(a,c)$ we obtain
\begin{equation}\begin{aligned}\label{Da:b}
0&=\sum_{k=1}^{m} R
\divided{e}{a-k+1}\KK{b_1+k-1}{b_2-k+1}\divided{f}{c-k}\\
+[a+1]&\divided{e}{a+1}\K \fc
+(-1)^m[a-b_1-m]\divided{e}{a-m}\KK{b_1+m}{b_2-m}\divided{f}{c-m-1}
\end{aligned}\end{equation}
where
\begin{align*}
R&=\vbinomm{b_1+k}{k}[a-k+1]-\vbinomm{b_1+k-1}{k-1}[a-b_1-k+1]\\ &
=[a+1]\vbinomm{b_1+k-1}{k}.
\end{align*}
Now the last term in \eqref{Da:b} is zero if $m=c$.
Otherwise $m=a<c$ and the term takes the form
\begin{align*} (-1)^a\vbinomm{b_1+a}{a}&[-b_1]\KK{b_1+a}{b_2-a}
\divided{f}{c-a-1}\\
&=
(-1)^{a+1}[a+1]\vbinomm{b_1+a}{a+1}\KK{b_1+a}{b_2-a}\divided{f}{c-a-1}.
\end{align*}
This shows that all the terms in equation \eqref{Da:b} have a
common factor of $[a+1]$. Putting these terms together and
dividing by $[a+1]$ we obtain the equality
$$0=\sum_{k=0}^M (-1)^k\vbinomm{b_1+k-1}{k}\divided{e}{a-k+1}
\KK{b_1+k-1}{b_2-k+1}\divided{f}{c-k}$$
where $M=m=\min(a+1,c)$ in the case $m=c$ and
$M=m+1=a+1=\min(a+1,c)$ otherwise. In either case $M=\min(a+1,c)$
and the induction is complete.
\end{pf}
\subsection{}\begin{xthm}\label{Db}
Suppose $a,b_1,b_2,c\in \N$ with $b_1+b_2=d$. Then if $s=a+b_1+c-d>0$
we have the equality
$$ \ea\K \fc
=\sum_{k=s}^{\min(a,c)}(-1)^{k-s}\vbinomm{k-1}{s-1}\vbinomm{b_1+k}{k}
\divided{e}{a-k}\KK{b_1+k}{b_2-k}\divided{f}{c-k}$$
and if
$s=a+b_2+c-d$ we have the equality
$$f^{(a)}\K e^{(c)}
= \sum_{k=s}^{\min(a,c)}(-1)^{k-s}\vbinomm{k-1}{s-1}\vbinomm{b_2+k}{k}
\divided{f}{a-k}\KK{b_1-k}{b_2+k}\divided{e}{c-k}.$$
\end{xthm}
\begin{pf}
By symmetry we need only verify the first equality.
We proceed by induction on $s$. The case $s=1$ is the content of
Theorem \ref{Da}. Let $a,b_1,b_2,c$ be given such that
$a+b_1+c-d=s+1$ and $b_1+b_2=d$. If $c<s$ then $a+b_1\geq d+1$
and so by Lemma \ref{Bf} we have
$$\ea\K\fc=0.$$
By a similar argument one sees that this also holds if $a<c$.
Hence we may assume that both $a$ and $c$ are $\geq s$. It is
enough to prove the result for the case $a\geq c\geq s$ since the
other case is similar.

Thus $a\geq 1$ and we have by induction
\begin{align*}
&\ea \K\fc= \frac{e}{[a]}\divided{e}{a-1}\K\fc\\
&= \frac{e}{[a]}\sum_{k=s}^{\min(a-1,c)}(-1)^{k-s}
\vbinomm{k-1}{s-1}\vbinomm{b_1+k}{k}
\divided{e}{a-1-k}\KK{b_1+k}{b_2-k}\divided{f}{c-k}\\
&=\frac{1}{[a]}\sum_{k=s}^{\min(a-1,c)}(-1)^{k-s}
\vbinomm{k-1}{s-1}\vbinomm{b_1+k}{k}
[a-k]\divided{e}{a-k}\KK{b_1+k}{b_2-k}\divided{f}{c-k}\\
&=\frac{[a-s]}{[a]}\vbinomm{b_1+s}{s}\divided{e}{a-s}\KK{b_1+s}{b_2-s}
\divided{f}{c-s}\\
&\ + \frac{1}{[a]}\sum_{k=s+1}^{\min(a-1,c)}(-1)^{k-s}
\vbinomm{k-1}{s-1}\vbinomm{b_1+k}{k}[a-k]
\divided{e}{a-k}\KK{b_1+k}{b_2-k}\divided{f}{c-k}.
\end{align*}
Now the first term of the last equality can be expanded by Theorem
\ref{Da} since $(a-s)+(b_1+s)+(c-s)=a+b_1+c-s=d+1$.
Putting this in and shifting the index of summation we obtain the
equalities
\begin{align*}
&e^{(a)} \K\fc\\
&=
\frac{[a-s]}{[a]}\vbinomm{b_1+s}{s}\sum_{k=1}^{c-s}(-1)^{k-1}\vbinomm{b_1+s+k}{k}
\divided{e}{a-s-k}\KK{b_1+s+k}{b_2-s+k}\divided{f}{c-s-k}\\
&\qquad -\frac{1}{[a]}\sum_{k=s+1}^{\min(a-1,c)}(-1)^{k-(s+1)}
\vbinomm{k-1}{s-1}\vbinomm{b_1+k}{k}[a-k]
\divided{e}{a-k}\KK{b_1+k}{b_2-k}\divided{f}{c-k}\\
&=
\frac{[a-s]}{[a]}\vbinomm{b_1+s}{s}\sum_{k=s+1}^{c}(-1)^{k-(s+1)}\vbinomm{b_1+k}{k-s}
\divided{e}{a-k}\KK{b_1+k}{b_2+k}\divided{f}{c-k}\\
&\qquad -\frac{1}{[a]}\sum_{k=s+1}^{\min(a-1,c)}(-1)^{k-(s+1)}
\vbinomm{k-1}{s-1}\vbinomm{b_1+k}{k}[a-k]
\divided{e}{a-k}\KK{b_1+k}{b_2-k}\divided{f}{c-k}
\end{align*}
Now the second term in the last equality above can be taken from
$s+1$ to $c$ since $\min(a-1,c)$ is different from $c$ only if
$a=c$, in which case the additional term in the sum will be zero
(the factor $[a-k]$ is zero when $k=c=a$). Putting the two sums
together and using the identity
$$
\frac{[a-s]}{[a]}\vbinomm{b_1+s}{s}\vbinomm{b_1+k}{k-s}-\frac{[a-k]}{[a]}
\vbinomm{k-1}{s-1}\vbinomm{b_1+k}{k}=\vbinomm{k-1}{s}\vbinomm{b_1+k}{k}.$$
we obtain
$$\ea\K\fc=\sum_{k=s+1}^c
(-1)^{k-(s+1)}\vbinomm{k-1}{s}\vbinomm{b_1+k}{k}\divided{e}{a-k}
\KK{b_1+k}{b_2-k}\divided{f}{c-k}$$
and this completes the induction.
\end{pf}

\subsection{}\label{Dc}
Theorem \ref{Db} combines with the isomorphism $\B_d\cong \Sv$ to
prove Theorem \ref{AAd}.

We now prove Theorem \ref{AAc}. First recall the representation
\begin{equation}\label{Dc:a}
\rho_d: \U\rightarrow \End_{\Q(v)}(E^{\ts d})\end{equation}
of section \ref{Ae}. Let $\Uz$ be the subalgebra of $\U$
generated by all $\divided{e}{m}$, $\divided{f}{m}$ ($m\in \N$)
and $\Ki^{\pm 1}$. The map in \eqref{Dc:a} restricts to give a map
$\Uz\rightarrow \End_{\Q(v)}(E^{\ts d})$. Let $E_{\A}$ be the
$\A$-submodule of $E$ spanned by the canonical basis elements
$e_1$ and $e_2$. It is clear that $E_\A$ is stable under the
action of $\Uz$ and hence so is $E_{\A}^{\ts d}$. Thus the image
of the above map is contain in $\End_{\A}(E_{\A}^{\ts d})$.
Consequently we have a representation
\begin{equation}\label{Dc:b}
\rho_d^{\A}: \Uz\rightarrow \End_{\A}(E_{\A}^{\ts d})\end{equation}
and by a result of Du \cite{Du} (see also \cite{Gr}) the
algebra $\Sz$ is precisely the image of this representation.

A fundamental result in \cite{Lu} (applied to the ${\mathfrak
{gl}}_2$ case) is that
$\Uz$ is a free $\A$-module with basis
\begin{equation}\label{Dc:c}
\left\{ e^{(a)} \Ka^{\delta_1} \Kb ^{\delta_2}\K f^{(c)}\quad | \quad
a,b_1,b_2,c\in \N,\quad \delta_i\in \{0,1\}.\right\}\end{equation}
Hence $\Sz$ is spanned over $\A$ by the images of the basis
elements of (\ref{Dc:c}). By Lemma \ref{Bd} and Theorem \ref{Be} the elements
$\Ka^{\delta_1} \Kb ^{\delta_2}\K$ with $\delta_i\in
\{0,1\}$ and $b_1,b_2$ arbitrary are $\A$-linear combinations
of those $\K$ with $b_1+b_2=d$. Therefore $\Sz$ is spanned by all
$\ea \K \fc$ with $b_1+b_2=d$. But by Theorem \ref{Db} we see that
the set of all such terms satisfying the condition $a+b_1+c\leq d$
is a spanning set for $\Sz$ over $\A$. Being linearly independent
over $\Q(v)$, this set is also linearly independent over $\A$. By symmetry,
the set $\fa \K \ec$ with $b_1+b_2=d$ and $a+b_2+c\leq d$ is linearly
independent over $\A$. This completes the proof of Theorem \ref{AAc}.

It remains to prove Theorem \ref{AAe}. 
By Lemma \ref{Bd} and Theorem \ref{Be} it follows that
\begin{align*}
\vbinomm{K_1}{b} = \vbinomm{K_1}{b}\cdot 1 
&= \vbinomm{K_1}{b}\sum_{b_1+b_2=d} \K\\
&= \sum_{b_1+b_2=d}\vbinomm{b_1}{b} \K
\end{align*}
and similarly we have
$$
\vbinomm{K_2}{b} = \sum_{b_1+b_2=d}\vbinomm{b_2}{b} \K.
$$
Moreover, the matrix of coefficients in these equations is, with
respect to an appropriate ordering of its rows and columns,
triangular with $1$'s on the main diagonal.  So these equations can be
inverted over $\A$ to obtain formulas expressing each $\K$
($b_1+b_2=d$) as an $\A$-linear combination of the $\vbinomm{K_1}{b}$
or the $\vbinomm{K_2}{b}$ as $b$ ranges from $0$ to $d$.  Thus it
follows that the sets $\{1, \vbinomm{K_1}{1}, \dots,
\vbinomm{K_1}{d}\}$ and $\{1, \vbinomm{K_2}{1}, \dots,
\vbinomm{K_2}{d}\}$ are spanning sets (over $\A$) for the algebra 
$S_\A^0(n,d)$.  Since it is already known that the rank of this free
$\A$-module is $d+1$ and since the ring $\A$ is commutative, it
follows that the two sets are in fact $\A$-bases for $S_\A^0(n,d)$.

\subsection{}\begin{xlem}\label{Dd}
Let $a,b,c$ be nonnegative integers satisfying $a+b+c>d$. Then the
elements $\ea \vbinomm{K_1}{b} \fc$ and $\fa \vbinomm{K_2}{b} \ec$ 
of degree $a+b+c$ are each expressible as $\A$-linear combinations 
of elements of the same form but of degree not exceeding $d$.
\end{xlem}

\begin{pf}
From the remarks preceding the lemma we know that $\vbinomm{K_1}{b}$
and $\vbinomm{K_2}{b}$ are expressible as $\A$-linear combinations of
the idempotents $\K$ ($b_1+b_2=d$).  Thus the elements $\ea
\vbinomm{K_1}{b} \fc$ (resp., $\fa \vbinomm{K_2}{b} \ec$) are expressible 
as $\A$-linear combinations of elements of the form
\begin{equation} \label{Dd:a}
\ea \K \fc  \quad (\text{resp., } \fa \K \ec)
\end{equation}
where $b_1+b_2=d$ and where the fake degree $a+b_1+c$ (resp.,
$a+b_2+c$) is strictly greater than $d$.  By Theorem \ref{AAd} it
follows that each element of the form \eqref{Dd:a} above is
expressible as an $\A$-linear combination of elements of
the same form
\begin{equation} \label{Dd:b}
\divided{e}{a'} K_{u_1,u_2} \divided{f}{c'}  \quad (\text{resp., } 
\divided{f}{a'} K_{u_1,u_2} \divided{e}{c'})
\end{equation}
where $a'<a$, $c'<c$, $u_1+u_2=d$, and $a'+u_1+c'\le d$ (resp.,
$a'+u_2+c'\le d$).  Now by expressing $K_{u_1,u_2}$ as an $\A$-linear
combination of elements of the form $\vbinomm{K_1}{b'}$, (resp.,
$\vbinomm{K_2}{b'}$) where $0\le b' \le d$ we obtain (via left and
right multiplication by appropriate elements) formulas expressing each
of the elements of the form \eqref{Dd:b} in terms of $\A$-linear
combinations of elements of the form
\begin{equation} \label{Dd:c}
\divided{e}{a'} \vbinomm{K_1}{b'} \divided{f}{c'}  \quad (\text{resp., } 
\divided{f}{a'} \vbinomm{K_2}{b'} \divided{e}{c'}).
\end{equation}
If an element of this form satisfies the constraint $a'+b'+c' \le d$
then we leave it be, but for those elements which do not satisfy this
constraint we repeat the entire process given above, replacing the
element by an $\A$-linear combination of elements of the same form, in
which for each element the degree of $e$ is strictly smaller, as is
the degree of $f$. After repeating the process finitely many times we
obtain the desired result.
\end{pf}

Now we can prove Theorem \ref{AAe}.  By Theorem \ref{AAc} we know that
the set
$$
B=\left\{\ea \K \fc \mid a+b_1+c\leq d,\ b_1+b_2=d\right\}
$$ 
is an $\A$-basis for $\Sz$. Now consider the set
$$
B'=\left\{\ea \vbinomm{K_1}{b} \fc \mid a+b+c\leq d\right\}.
$$ 
The sets $B$ and $B'$ have the same cardinality and the ring $\A$ is
commutative, thus it will follow that $B'$ is a basis once we can show
that it spans.

We know that $\Sz$ is spanned by elements of the form $\ea
\vbinomm{K_1}{b} \fc$ since $S_\A = S_\A^+ S_\A^0 S_\A^-$.  By the
above lemma we know that each such element not satisfying the
constraint $a+b+c\le d$ is expressible as an $\A$-linear combination
of elements which do satisfy that constraint.  It follows that the set
$B'$ is a spanning set for $\Sz$. This proves the first part of
Theorem \ref{AAe}.  The second part of the theorem follows by
symmetry.

%\newpage
\bigskip

\parskip=0pt\parindent=0pt\sf
\par Mathematical and Computer Sciences
\par Loyola University Chicago
\par Chicago, Illinois 60626 U.S.A.
\par E-mail: doty@math.luc.edu
\par \ \ \ \ \ \ \ \ \ tonyg@math.luc.edu

\end{document}